\input amstex
\documentstyle{amsppt}
%----------------------------------------------------------------
% Title:     On deformations of metrics and their associated 
%            spinor structures.
% Author:    Ruslan Sharipov
% Comments:  AmSTeX, 22 pages, amsppt style
% MSC-class: 53B30, 81T20, 83F05
%----------------------------------------------------------------
%           Replacement for output macro definition
%
\catcode`@=11
\redefine\output@{%
  \def\break{\penalty-\@M}\let\par\endgraf
  \ifodd\pageno\global\hoffset=105pt\else\global\hoffset=8pt\fi  
  \shipout\vbox{%
    \ifplain@
      \let\makeheadline\relax \let\makefootline\relax
    \else
      \iffirstpage@ \global\firstpage@false
        \let\rightheadline\frheadline
        \let\leftheadline\flheadline
      \else
        \ifrunheads@ %\let\makefootline\relax
        \else \let\makeheadline\relax
        \fi
      \fi
    \fi
    \makeheadline \pagebody \makefootline}%
  \advancepageno \ifnum\outputpenalty>-\@MM\else\dosupereject\fi
}
\catcode`\@=\active
%----------------------------------------------------------------
\nopagenumbers
\def\negskp{\hskip -2pt}
\def\Alpha{\operatorname{A}}
\def\MatGrSO{\operatorname{SO}}
\def\MatGrSL{\operatorname{SL}}
\def\MatGrGL{\operatorname{GL}}
\def\const{\operatorname{const}}
\def\idop{\operatorname{\bold{id}}}
\accentedsymbol\hatboldd{\kern 2pt\hat{\kern -2pt\bold d}}
\accentedsymbol\hatboldD{\kern -1pt\hat{\kern 1pt\bold D}}
\def\vtrule{\vrule height 12pt depth 6pt}
\def\vtttrule{\vrule height 12pt depth 19pt}
\def\boxit#1#2{\vcenter{\hsize=122pt\offinterlineskip\hrule
  \line{\vtttrule\hss\vtop{\hsize=120pt\centerline{#1}\vskip 5pt
  \centerline{#2}}\hss\vtttrule}\hrule}}
\def\blue#1{#1}
\catcode`#=11\def\diez{#}\catcode`#=6
\catcode`_=11\def\podcherkivanie{_}\catcode`_=8
%\catcode`~=11\def\volna{~}\catcode`~=\active
\def\mycite#1{\cite{\blue{#1}}\immediate\special{ps:
     ShrHPSdict begin /ShrBORDERthickness 0 def}}

\def\mytag#1{%
    \tag#1}
\def\mythetag#1{\thetag{\blue{#1}}\immediate\special{ps:
     ShrHPSdict begin /ShrBORDERthickness 0 def}}
\def\myrefno#1{\no#1}
\def\myhref#1#2{\blue{#2}\immediate\special{ps:
     ShrHPSdict begin /ShrBORDERthickness 0 def}}
\def\myEarXivlink{\myhref{http://arXiv.org}{http:/\negskp/arXiv.org}}
\def\myGeoCities{\myhref{http://www.geocities.com}{GeoCities}}
\def\mytheorem#1{\csname proclaim\endcsname{Theorem #1}}
\def\mythetheorem#1{\blue{#1}\immediate\special{ps:
     ShrHPSdict begin /ShrBORDERthickness 0 def}}
\def\mylemma#1{\csname proclaim\endcsname{Lemma #1}}

\def\mycorollary#1{\csname proclaim\endcsname{Corollary #1}}

\def\mydefinition#1{\definition{Definition #1}}
\def\mythedefinition#1{\blue{#1}\immediate\special{ps:
     ShrHPSdict begin /ShrBORDERthickness 0 def}}

%----------------------------------------------------------------
% Cyrillic fonts definition
%\font\eightcyr=wncyr8
%----------------------------------------------------------------
\pagewidth{360pt}
\pageheight{606pt}
\topmatter
\title
On deformations of metrics and their associated 
spinor structures.
\endtitle
\author
R.~A.~Sharipov
\endauthor
\address 5 Rabochaya street, 450003 Ufa, Russia\newline
\vphantom{a}\kern 12pt Cell Phone: +7(917)476 93 48
\endaddress
\email \vtop to 30pt{\hsize=280pt\noindent
\myhref{mailto:r-sharipov\@mail.ru}
{r-sharipov\@mail.ru}\newline
\myhref{mailto:R\podcherkivanie Sharipov\@ic.bashedu.ru}
{R\_\hskip 1pt Sharipov\@ic.bashedu.ru}\vss}
\endemail
\urladdr
\vtop to 20pt{\hsize=280pt\noindent
\myhref{http://www.geocities.com/r-sharipov}
{http:/\negskp/www.geocities.com/r-sharipov}\newline
\myhref{http://www.freetextbooks.boom.ru/index.html}
{http:/\negskp/www.freetextbooks.boom.ru/index.html}\vss}
\endurladdr
\abstract
    Smooth deformations of a Minkowski type metric in a 
four-dimensional space-time manifold are considered. 
Deformations of the basic spin-tensorial fields associated
with this metric are calculated and their application to 
calculating the energy-momentum tensor of a massive spin
1/2 particle is shown.
\endabstract
\endtopmatter
\loadbold
\loadeufb
\TagsOnRight
\document
\accentedsymbol\hatboldsymbolgamma{\kern 1.3pt\hat{\kern -1.3pt
   \boldsymbol\gamma}}
\accentedsymbol\hatgamma{\kern 1.3pt\hat{\kern -1.3pt\gamma}}

\rightheadtext{On deformations of metrics \dots}
\head
1. Introduction. 
\endhead
    Let $M$ be a four-dimensional orientable manifold equipped 
with a Minkowski type metric $\bold g$ and with a 
polarization\footnote{ \ A polarization is a discrete geometric 
structure (like an orientation) that marks the future half light 
cone at each point of $M$ (see more details in \mycite{1}).}.
\adjustfootnotemark{-1}
In general relativity such a manifold $M$ is used as a stage for
all physical phenomena. When describing the spin phenomenon $M$
is additionally assumed to be a spin manifold. In this case it
admits two spinor bundles: the bundle of Weyl spinors $SM$ and
the bundle of Dirac spinors $DM$ (see definitions below).\par
     Relativistic quantum particles in $M$ are described by their
fields, while fields are introduced through their contribution
to the action 
$$
\hskip -2em
S=S_{\text{gravity}}+S_{\text{matter}}.
\mytag{1.1}
$$
The metric $\bold g$ itself is interpreted as a gravitation field, 
while other fields are called matter fields. The influence of the 
matter upon the gravitation field is described by the following 
Einstein equation:
$$
R_{ij}-\frac{R}{2}\,g_{ij}=\frac{8\,\pi\,G}{c^4}T_{ij}.
\mytag{1.2}
$$
Here $R_{ij}$ are the components of the Ricci tensor and $R$ is the
scalar curvature (both are produced from the components of the metric
tensor $g_{ij}$ and their derivatives). The right hand side of
\mythetag{1.2} comprises two fundamental physical constants: the 
Newtonian constant of gravitation $G$ and the light velocity in vacuum
$c$. Here are their numeric values taken from the 
\myhref{http://physics.nist.gov/cuu/Constants}{NIST} site:
$$
\xalignat 2
&G\approx 6.67428\cdot 10^{-8}\ \text{\it cm}^3\ \text{\it\! g}^{-1}
\ \text{\it\! sec}^{-2},
&&c\approx 2.99792458\cdot 10^{10}\ \text{\it cm}
\ \ \text{\it\! sec}^{-1}.
\qquad
\mytag{1.3}
\endxalignat
$$
Apart from the constants \mythetag{1.3}, the right hand side of 
the equation \mythetag{1.2} comprises the components of the
energy-momentum tensor $T_{ij}$. They are defined through the 
variational derivatives of the matter action in \mythetag{1.1}:
$$
\hskip -2em
T_{ij}=2\,c\,\frac{\delta S_{\text{matter}}}{\delta g^{ij}}.
\mytag{1.4}
$$
The main goal of this paper is to clarify the procedure of 
calculating the right hand side of the formula \mythetag{1.4}
in the case where matter fields have spinor components.
\head
2. The bundle of Weyl spinors.
\endhead
     The bundle of Weyl spinors $SM$ is a special two-dimensional 
complex vector bundle over the space-time manifold $M$. It is a 
special bundle since it is related in a special way to the tangent
bundle $TM$. This relation is based on the well-known algebraic
fact --- the group homomorphism
$$
\hskip -2em
\phi:\ \MatGrSL(2,\Bbb C)\to\MatGrSO^+(1,3,\Bbb R).
\mytag{2.1}
$$
The group $\MatGrSO^+(1,3,\Bbb R)$ in the right hand side of 
\mythetag{2.1} is the special orthochronous Lorentz group. It 
is linked to the metric structure of the manifold $M$ through
frames. In the context of general relativity they are also
called {\bf vierbeins}. This German word means ``four feet'' 
expressing the idea of frames in the case of a four-dimensional
space-time manifold.
\mydefinition{2.1} A frame of the tangent bundle $TM$ for the 
space-time manifold $M$ is an ordered quadruple of smooth 
vector fields $\boldsymbol\Upsilon_0,\,\boldsymbol\Upsilon_1,
\,\boldsymbol\Upsilon_2,\,\boldsymbol\Upsilon_3$ defined in
some open domain $U\subset M$ and linearly independent at each
point $p\in U$.
\enddefinition
     For any two frames $\boldsymbol\Upsilon_0,\,\boldsymbol
\Upsilon_1,\,\boldsymbol\Upsilon_2,\,\boldsymbol\Upsilon_3$ and
$\tilde{\boldsymbol\Upsilon}_0,\,\tilde{\boldsymbol\Upsilon}_1,
\,\tilde{\boldsymbol\Upsilon}_2,\,\tilde{\boldsymbol\Upsilon}_3$
with intersecting domains the transition matrices arise. At each 
point $p\in U\cap\tilde U$ we have
$$
\xalignat 2
&\hskip -2em
\tilde{\boldsymbol\Upsilon}_i=\sum^3_{j=0}S^j_i\,\boldsymbol\Upsilon_j,
&&\boldsymbol\Upsilon_i=\sum^3_{j=0}T^j_i\,\tilde{\boldsymbol\Upsilon}_j.
\mytag{2.2}
\endxalignat
$$
\mydefinition{2.2} A frame $\boldsymbol\Upsilon_0,\,\boldsymbol
\Upsilon_1,\,\boldsymbol\Upsilon_2,\,\boldsymbol\Upsilon_3$ of the
tangent bundle $TM$ is called a holonomic frame if its vector fields
$\boldsymbol\Upsilon_0,\,\boldsymbol\Upsilon_1,\,\boldsymbol\Upsilon_2,
\,\boldsymbol\Upsilon_3$ commute with each other. Otherwise it is called 
a non-holonomic frame:
$$
\hskip -2em
[\boldsymbol\Upsilon_i,\boldsymbol\Upsilon_j]=\sum^3_{k=0}
c^{\,k}_{ij}\,\boldsymbol\Upsilon_k.
\mytag{2.3}
$$
The coefficients $c^{\,k}_{ij}$ in \mythetag{2.3} are called 
the commutations coefficients of this non-holonomic frame.
\enddefinition
    Note that each local chart with coordinates $x^0,\,x^1,\,x^2,
\,x^3$ in $M$ produces its associated holonomic frame of the 
coordinate vector fields: 
$$
\xalignat 4
&\hskip -2em
\boldsymbol\Upsilon_0=\frac{\partial}{\partial x^0},
&&\boldsymbol\Upsilon_1=\frac{\partial}{\partial x^1},
&&\boldsymbol\Upsilon_2=\frac{\partial}{\partial x^2},
&&\boldsymbol\Upsilon_3=\frac{\partial}{\partial x^3}.
\qquad
\mytag{2.4}
\endxalignat
$$
And conversely, each holonomic frame can be represented in the form
\pagebreak of \mythetag{2.4} in some neighborhood of each point where 
it is defined. However, in general case frames are non-holonomic.
\par
    Frames are used to represent vectorial and tensorial fields 
in the coordinate form. In particular, if we take the metric tensor
$\bold g$, it is represented by a square matrix $g_{ij}$. Its dual
metric tensor is represented by the inverse matrix $g^{ij}$.
\mydefinition{2.3} A frame $\boldsymbol\Upsilon_0,\,\boldsymbol
\Upsilon_1,\,\boldsymbol\Upsilon_2,\,\boldsymbol\Upsilon_3$ of the
tangent bundle $TM$ is called an orthonormal frame if the components
of the metric tensor $\bold g$ and its dual metric tensor are given 
by the standard Minkowski matrix 
$$
\hskip -2em
g_{ij}=g^{ij}=\Vmatrix 1 & 0 & 0 & 0\\0 & -1 & 0 & 0\\
0 & 0 & -1 & 0\\0 & 0 & 0 & -1\endVmatrix.
\mytag{2.5}
$$
\enddefinition
     Remember that we assume $M$ to be an orientable manifold. 
Therefore, we can fix an orientation on it. Then all frames will
subdivide into two types --- right frames and left frames. If both
frames in \mythetag{2.2} are right or if both of them are left,
then we say that they have the same orientation. In this case
$$
\xalignat 2
&\det S>0, &&\det T>0.
\mytag{2.6}
\endxalignat 
$$
Otherwise we say that these frames have opposite orientations.
In this case 
$$
\xalignat 2
&\det S<0, &&\det T<0.
\mytag{2.7}
\endxalignat 
$$
If, additionally, both frames in \mythetag{2.2} are orthonormal, then
$S$ and $T$ are Lorentzian matrices. For their determinants we have
$$
\det S=\det T=\pm 1.
\mytag{2.8}
$$
The sign in \mythetag{2.8} is chosen according to \mythetag{2.6} or
according to \mythetag{2.7} depending on the orientations of frames.
\mydefinition{2.4} A frame $\boldsymbol\Upsilon_0,\,\boldsymbol
\Upsilon_1,\,\boldsymbol\Upsilon_2,\,\boldsymbol\Upsilon_3$ of the
tangent bundle $TM$ is called a positively polarized frame if its
first vector $\boldsymbol\Upsilon_0$ lies inside the future half
light cone, i\.\,e\. if it is a time-like vector directed to the
future.
\enddefinition
     Assume that both frames in \mythetag{2.2} are positively 
polarized right orthonormal frames. In this case both transition
matrices $S$ and $T$ in \mythetag{2.2} are Lorentzian matrices 
with the unit determinant such that $S^0_0>0$ and $T^0_0>0$.
This means that
$$
\xalignat 2
&\hskip -2em
S\in\MatGrSO^+(1,3,\Bbb R),
&&T\in\MatGrSO^+(1,3,\Bbb R).
\mytag{2.9}
\endxalignat
$$
Thus, having postulated the presence three geometric structures in 
$M$ --- the metric, the orientation, and the polarization, we have
implemented the special orthochronous Lorentz group through the
frame transition matrices \mythetag{2.9}. This means that the 
structural group of the tangent bundle $TM$ reduces from 
$\MatGrGL(4,\Bbb R)$ to $\MatGrSO^+(1,3,\Bbb R)$. Generally
speaking, \pagebreak any reduction of the structural group of $TM$ 
from $\MatGrGL(4,\Bbb R)$ to some subgroup of $\MatGrGL(4,\Bbb R)$ 
is due to the presence of some definite geometric structures in $M$. 
Some similar fact can be formulated for any other bundle over the 
base manifold $M$.\par
\mydefinition{2.5} A frame of the spinor bundle $SM$ for the 
space-time manifold $M$ is an ordered pair of smooth 
spinor fields $\boldsymbol\Psi_1,\,\boldsymbol\Psi_2$ defined in
some open domain $U\subset M$ and linearly independent at each
point $p\in U$.
\enddefinition
     For any two frames $\boldsymbol\Psi_1,\,\boldsymbol
\Psi_2$ and $\tilde{\boldsymbol\Psi}_1,\,\tilde{\boldsymbol\Psi}_2$
with intersecting domains the transition matrices arise. At each 
point $p\in U\cap\tilde U$ we have
$$
\xalignat 2
&\hskip -2em
\tilde{\boldsymbol\Psi}_i=\sum^2_{j=1}\goth S^j_i\,\boldsymbol\Psi_j,
&&\boldsymbol\Psi_i=\sum^2_{j=1}\goth T^j_i\,\tilde{\boldsymbol\Psi}_j.
\mytag{2.10}
\endxalignat
$$
\mydefinition{2.6} A two-dimensional complex vector bundle $SM$
over the space-time manifold $M$ is called the bundle of Weyl
spinors if for any positively polarized right orthonormal frame
$\boldsymbol\Upsilon_0,\,\boldsymbol\Upsilon_1,\,\boldsymbol
\Upsilon_2,\,\boldsymbol\Upsilon_3$ and for any point $p$ of its
domain $U$ some frame $\boldsymbol\Psi_1,\,\boldsymbol\Psi_2$ of 
$SM$ defined in $U$ or in some smaller neighborhood of the point $p$
is canonically associated to $\boldsymbol\Upsilon_0,\,\boldsymbol
\Upsilon_1,\,\boldsymbol\Upsilon_2,\,\boldsymbol\Upsilon_3$ in such 
a way that for any two positively polarized right orthonormal frames
of $TM$ related by the matrices $S$ and $T$ in \mythetag{2.2} their
associated frames in $SM$ are related to each other by the matrices
$\goth S\in\MatGrSL(2,\Bbb C)$ and $\goth T\in\MatGrSL(2,\Bbb C)$ 
in \mythetag{2.10}, where $S=\phi(\goth S)$ and $T=\phi(\goth T)$ and 
where $\phi$ is the group homomorphism \mythetag{2.1}.
\enddefinition
     According to the definition~\mythedefinition{2.6} the structural
group of the spinor bundle $SM$ reduces from $\MatGrGL(2,\Bbb C)$ to
$\MatGrSL(2,\Bbb C)$. Due to this reduction the bundle of Weyl spinors
is equipped with the skew-symmetric spinor metric $\bold d$. It is
given by the matrix 
$$
\hskip -2em
d_{ij}=\Vmatrix 0 & 1\\ 
\vspace{1ex} -1 & 0\endVmatrix
\mytag{2.11}
$$
in any frame $\boldsymbol\Psi_1,\,\boldsymbol\Psi_2$ canonically
associated with some positively polarized right orthonormal frame 
of $TM$. The dual metric for \mythetag{2.11} in such a frame is 
given by the matrix inverse to the matrix \mythetag{2.11}:
$$
\hskip -2em
d^{ij}=\Vmatrix 0 & -1\\ 
\vspace{1ex} 1 & 0\endVmatrix.
\mytag{2.12}
$$
\mydefinition{2.7} A frame $\boldsymbol\Psi_1,\,\boldsymbol\Psi_2$
of the spinor bundle $SM$ is called an orthonormal frame if the
spinor metric $\bold d$ and its dual metric are given by the matrices 
\mythetag{2.11} and \mythetag{2.12} in this frame.
\enddefinition
    As we see, any frame $\boldsymbol\Psi_1,\,\boldsymbol\Psi_2$
canonically associated with some positively polarized right orthonormal
frame $\boldsymbol\Upsilon_0,\,\boldsymbol\Upsilon_1,\,\boldsymbol
\Upsilon_2,\,\boldsymbol\Upsilon_3$ of $TM$ is an orthonormal frame
in the sense of the definition~\mythedefinition{2.7}. The converse 
proposition is also true, i\.\,e\. any orthonormal frame $\boldsymbol
\Psi_1,\,\boldsymbol\Psi_2$ of $SM$ is canonically associated with 
some positively polarized right orthonormal frame in $TM$. Therefore, 
we use the following diagram:
$$
\hskip -2em
\aligned
&\boxit{\lower 5pt\hbox{Orthonormal frames}}{}\to
\boxit{Positively polarized}{right orthonormal frames}
\endaligned\quad
\mytag{2.13}
$$\par
     The reduction of $\MatGrGL(2,\Bbb C)$ to $\MatGrSL(2,\Bbb C)$ for
Weyl spinors is concordant with the reduction of $\MatGrGL(4,\Bbb R)$ 
to $\MatGrSO^+(1,3,\Bbb R)$ in $TM$. For this reason the spinor 
bundle $SM$ has one more geometric structure given by the Infeld-van
der Waerden field $\bold G$. For any frame pair $\boldsymbol\Upsilon_0,
\,\boldsymbol\Upsilon_1,\,\boldsymbol\Upsilon_2,\,\boldsymbol\Upsilon_3$
and $\boldsymbol\Psi_1,\,\boldsymbol\Psi_2$ canonically associated to
each other according to the diagram \mythetag{2.13} the components of 
the Infeld-van der Waerden field are given by the following four 
Pauli matrices:
$$
\xalignat 2
&\hskip -2em
G^{i\kern 0.5pt\bar i}_0=\Vmatrix 1 & 0\\ 0 & 1\endVmatrix
=\sigma_0,
&&G^{i\kern 0.5pt\bar i}_2=\Vmatrix 0 & -i\\ i & 0\endVmatrix
=\sigma_2,\\
\vspace{-1.4ex}
&&&\mytag{2.14}\\
\vspace{-1.4ex}
&\hskip -2em
G^{i\kern 0.5pt\bar i}_1=\Vmatrix 0 & 1\\ 1 & 0\endVmatrix
=\sigma_1,
&&G^{i\kern 0.5pt\bar i}_3=\Vmatrix 1 & 0\\ 0 & -1\endVmatrix
=\sigma_3.
\endxalignat
$$
The barred letter $\bar i$ in \mythetag{2.14} is used in order to
indicate that the second upper index is a conjugate spinor index.
The first upper index $i$ in \mythetag{2.14} is a regular spinor 
index, while the lower index is a spacial one. It runs from $0$
to $3$ and enumerates the Pauli matrices \mythetag{2.14}. The 
following table summarize the spin tensorial types of the basic
fields $\bold d$ and $\bold G$ of the bundle of Weyl spinors:
$$
\vcenter{\hsize 10cm
\offinterlineskip\settabs\+\indent
\vtrule
\hskip 1.2cm &\vtrule % Quantity
\hskip 5.2cm &\vtrule % Unit
\hskip 2.8cm &\vtrule % Relation
\cr\hrule 
\+\vtrule
\hfill\,Symbol\hfill&\vtrule
\hfill Name\hfill &\vtrule
\hfill Spin-tensorial\hfill &\vtrule\cr
\vskip -0.2cm
\+\vtrule
\hfill &\vtrule
\hfill \hfill&\vtrule
\hfill type\hfill&\vtrule\cr\hrule
\+\vtrule
\hfill $\bold d$\hfill&\vtrule
\hfill Skew-symmetric metric tensor\hfill&\vtrule
\hfill $(0,2|0,0|0,0)$\hfill&\vtrule\cr\hrule
\+\vtrule
\hfill$\bold G$\hfill&\vtrule
\hfill Infeld-van der Waerden field\hfill&\vtrule
\hfill $(1,0|1,0|0,1)$\hfill&\vtrule\cr\hrule
}\quad
\mytag{2.15}
$$
The third column of the table \mythetag{2.15} says that in the
coordinate form the metric tensor $\bold d$ has two lower spinor 
indices, while the Infeld-van der Waerden field $\bold G$ has one 
upper spinor index, one upper conjugate spinor index, and one lower
spacial index. These facts are in agreement with \mythetag{2.11} 
and with \mythetag{2.14}.
\head
3. Deformations of the metric.
\endhead
     The space-time manifold $M$ is equipped with the metric $\bold g$.
Let's change this metric replacing it by some other metric $\hat\bold g$
with the same signature. Locally, in the coordinate form, at each fixed
point $p\in M$ the metrics $\bold g$ and $\hat\bold g$ are related to 
each other by means of the following formula:
$$
\hskip -2em
g_{ij}=\sum^3_{p=0}\sum^3_{q=0}F^p_i\,F^q_j\,\hat g_{p\kern 0.3pt q}.
\mytag{3.1}
$$
This formula reflects the purely algebraic fact that each quadratic 
form can be diagonalized (see \mycite{2}). We need smooth deformations
of the metric $\bold g$, i\.\,e\. smoothly depending on a point $p\in 
M$. For this reason we assume $F^p_i$ and $F^q_j$ in \mythetag{3.1}
to be the components of some smooth non-degenerate operator-valued field 
$\bold F$. We need small deformations of metric in \mythetag{1.4}. 
Therefore we assume $\bold F$ to be close to the identity operator 
$\idop$. This fact is expressed by the formula
$$
\hskip -2em
\bold F=\idop+\,\varepsilon\,\bold f.
\mytag{3.2}
$$
Here $\varepsilon\to 0$ is a small parameter. For $F^p_i$ and $F^q_j$
the formula \mythetag{3.2} yields
$$
\xalignat 2
&\hskip -2em
F^p_i=\delta^{\,p}_i+\varepsilon\,f^p_i,
&&F^q_j=\delta^{\,q}_j+\varepsilon\,f^q_j.
\mytag{3.3}
\endxalignat 
$$
Substituting \mythetag{3.3} into the formula \mythetag{3.1}, we obtain
$$
\hskip -2em
g_{ij}=\hat g_{ij}+\sum^3_{p=0}\varepsilon\,f^p_i\,\hat g_{pj}
+\sum^3_{q=0}\varepsilon\,\hat g_{iq}\,f^q_j+\dots\,.
\mytag{3.4}
$$
By dots in \mythetag{3.4} we denote the higher order terms with respect 
to the small parameter $\varepsilon$. We see that $g_{ij}\sim\hat g_{ij}$
as $\varepsilon\to 0$. For this reason we can set $\hat g_{pj}=g_{pj}$
and $\hat g_{iq}=g_{iq}$ in \mythetag{3.4}. As a result the formula
\mythetag{3.4} is written as 
$$
\hskip -2em
\hat g_{ij}=g_{ij}-\sum^3_{p=0}\varepsilon\,f^p_i\,g_{pj}
-\sum^3_{q=0}\varepsilon\,g_{iq}\,f^q_j+\dots\,.
\mytag{3.5}
$$
In order to simplify the formula \mythetag{3.5} we denote
$$
\hskip -2em
f_{ij}=\sum^3_{q=0}g_{iq}\,f^q_j.
\mytag{3.6}
$$
The formula \mythetag{3.6} is a particular instance of the standard 
index lowering procedure. It produces the twice covariant tensor field
with the components $f_{ij}$ from the operator-valued field $\bold f$ 
used in \mythetag{3.2}. Applying \mythetag{3.6} to \mythetag{3.5} we 
get
$$
\hskip -2em
\hat g_{ij}=g_{ij}-\varepsilon\,(f_{ij}+f_{j\kern 0.5pt i})+\dots\,.
\mytag{3.7}
$$
The formula \mythetag{3.7} is equivalent to the following one:
$$
\hskip -2em
\delta g_{ij}=-\varepsilon\,(f_{ij}+f_{j\kern 0.5pt i}).
\mytag{3.8}
$$
Thus, we see that only the $\bold g$-symmetric part of the operator
field $\bold f$ is actually contribute to the variation of the metric 
$\bold g$ under the transformation \mythetag{3.1}. For this reason
we choose $\bold g$-symmetric operator $\bold f$ by setting
$$
\hskip -2em
f^i_j=\frac{1}{2}\sum^3_{p=0}g^{ip}\,h_{pj}
\mytag{3.9}
$$
and assuming $h_{ij}$ to be symmetric:
$$
\hskip -2em
h_{ij}=h_{j\kern 0.5pt i}.
\mytag{3.10}
$$
Then from \mythetag{3.8}, \mythetag{3.9}, and \mythetag{3.10} we derive
$$
\hskip -2em
\delta g_{ij}=-\varepsilon\,h_{ij}.
\mytag{3.11}
$$\par
    Now let's consider the components of the dual metric tensor
$\hat g^{ij}$. They form the matrix inverse to $\hat g_{ij}$. For 
the component of this inverse matrix from \mythetag{3.5} we derive
$$
\hskip -2em
\hat g^{ij}=g^{ij}+\sum^3_{q=0}\varepsilon\,f^i_q\,g^{qj}
+\sum^3_{p=0}\varepsilon\,g^{ip}\,f^j_p+\dots\,.
\mytag{3.12}
$$
Like in \mythetag{3.6}, for the sake of brevity we denote
$$
\hskip -2em
f^{ij}=\sum^3_{q=0}f^i_q\,g^{qj}.
\mytag{3.13}
$$
Then, substituting \mythetag{3.13} into \mythetag{3.12}, we transform
\mythetag{3.12} as follows:
$$
\hskip -2em
\hat g^{ij}=g^{ij}+\varepsilon\,(f^{ij}+f^{j\kern 0.5pt i})
+\dots\,.
\mytag{3.14}
$$
Applying \mythetag{3.9} to \mythetag{3.13} and taking into account
the symmetry $g^{qj}=g^{j\kern 0.3pt q}$, we get
$$
\hskip -2em
f^{ij}=\frac{1}{2}\sum^3_{p=0}\sum^3_{q=0}g^{i\kern 0.5pt p}\,g^{jq}
\,h_{p\kern 0.3pt q}.
\mytag{3.15}
$$
Relying on \mythetag{3.14} and \mythetag{3.15}, it is convenient 
to denote
$$
\hskip -2em
h^{ij}=\sum^3_{p=0}\sum^3_{q=0}g^{i\kern 0.5pt p}\,g^{jq}
\,h_{p\kern 0.3pt q}.
\mytag{3.16}
$$
In terms of \mythetag{3.16}, the formula \mythetag{3.14} is written as
$$
\hskip -2em
\hat g^{ij}=g^{ij}+\varepsilon\,h^{ij}+\dots\,.
\mytag{3.17}
$$
Due to \mythetag{3.17} the variations of the components $g^{ij}$ 
of the dual metric tensor under the transformation \mythetag{3.1} 
are given by the formula
$$
\hskip -2em
\delta g^{ij}=\varepsilon\,h^{ij},
\mytag{3.18}
$$
where $h^{ij}$ are produced from $h_{ij}$ used in \mythetag{3.11}
by means of the standard index raising procedure \mythetag{3.16}.
For this reason $h^{ij}$ in \mythetag{3.18} are symmetric:
$$
\hskip -2em
h^{ij}=h^{j\kern 0.5pt i}.
\mytag{3.19}
$$
The symmetry \mythetag{3.19} is concordant with the symmetry of the 
metric tensor $\bold g$ itself.
\head
4. Deformations of the spinor structure.
\endhead
    In our approach, the spinor bundle $SM$ is introduced through 
the orthonormal frames (see the definition~\mythedefinition{2.6}
above). Therefore it is related to the metric $\bold g$. Let's 
study this relation. Assume that $\boldsymbol\Upsilon_0,
\,\boldsymbol\Upsilon_1,\,\boldsymbol\Upsilon_2,\,\boldsymbol
\Upsilon_3$ is some positively polarized right orthonormal frame
of the metric $\bold g$. Then the matrix $g_{ij}$ is diagonal and 
coincides with \mythetag{2.5} in this frame. Using the operator
$\bold F$ we define the other frame $\hat{\boldsymbol\Upsilon}_0,
\,\hat{\boldsymbol\Upsilon}_1,\,\hat{\boldsymbol\Upsilon}_2,
\,\hat{\boldsymbol\Upsilon}_3$ introducing it by means of the formula
$$
\hskip -2em
\hat{\boldsymbol\Upsilon}_i=\sum^3_{j=0}F^j_i\,\boldsymbol\Upsilon_j.
\mytag{4.1}
$$
The components of the operator $\bold F$ in \mythetag{4.1} act just
like the components of the transition matrix $S$ in \mythetag{2.2}.
The formula \mythetag{4.1} is equivalent to
$$
\hskip -2em
\hat{\boldsymbol\Upsilon}_i=\bold F(\boldsymbol\Upsilon_i).
\mytag{4.2}
$$
Using \mythetag{3.1}, \mythetag{2.2} and combining them with 
\mythetag{4.1} or with \mythetag{4.2}, one easily proves that
$g_{ij}$ coincides with the matrix of the deformed metric 
$\hat\bold g$ represented in the frame $\hat{\boldsymbol\Upsilon}_0,
\,\hat{\boldsymbol\Upsilon}_1,\,\hat{\boldsymbol\Upsilon}_2,
\,\hat{\boldsymbol\Upsilon}_3$. Therefore the metric $\hat\bold g$
is diagonal in the frame $\hat{\boldsymbol\Upsilon}_0,
\,\hat{\boldsymbol\Upsilon}_1,\,\hat{\boldsymbol\Upsilon}_2,
\,\hat{\boldsymbol\Upsilon}_3$ and its matrix coincides with 
\mythetag{2.5} in this frame. In other words, the frame
$\hat{\boldsymbol\Upsilon}_0,\,\hat{\boldsymbol\Upsilon}_1,
\,\hat{\boldsymbol\Upsilon}_2,\,\hat{\boldsymbol\Upsilon}_3$
is an orthonormal frame for the deformed metric $\hat\bold g$.
Moreover, from \mythetag{3.2} and \mythetag{3.3} we derive that
$$
\xalignat 2
&\hskip -2em
\det\bold F\to 1,
&&F^0_0\to 1,
\mytag{4.3}
\endxalignat
$$
as $\varepsilon\to 0$. The relationships \mythetag{4.3} mean that
for sufficiently small values of $\varepsilon$ the determinant of
the operator $\bold F$ and its component $F^0_0$ both are positive:
$$
\xalignat 2
&\hskip -2em
\det\bold F>0,
&&F^0_0>0.
\mytag{4.4}
\endxalignat
$$
Due to \mythetag{4.4} the deformed frame $\hat{\boldsymbol\Upsilon}_0,
\,\hat{\boldsymbol\Upsilon}_1,
\,\hat{\boldsymbol\Upsilon}_2,\,\hat{\boldsymbol\Upsilon}_3$
is a positively polarized and right frame. Thus, using the 
deformation operator $\bold F$, we can produce a positively
polarized right orthonormal frame for the deformed metric
$\hat\bold g$ from any positively polarized right orthonormal 
frame of the initial metric $\bold g$.\par
     Let's recall that, according to the
definition~\mythedefinition{2.6}, each positively polarized right
orthonormal frame $\boldsymbol\Upsilon_0,\,\boldsymbol\Upsilon_1,
\,\boldsymbol\Upsilon_2,\,\boldsymbol\Upsilon_3$ is associated with
some frame $\boldsymbol\Psi_1,\,\boldsymbol\Psi_2$ of the spinor
bundle $SM$. Let's declare the deformed frame $\hat{\boldsymbol
\Upsilon}_0,\,\hat{\boldsymbol\Upsilon}_1,\,\hat{\boldsymbol
\Upsilon}_2,\,\hat{\boldsymbol\Upsilon}_3$ to be associated with
the same spinor frame $\boldsymbol\Psi_1,\,\boldsymbol\Psi_2$. As 
a result we get the diagram
$$
\vcenter to 56pt{\vss\includegraphics{MetDef01.eps}\vss}
\mytag{4.5}
$$
The arched arrow on the diagram \mythetag{4.5} corresponds to the 
frame deformation, while the horizontal arrows are frame associations, 
the red arrow being our newly declared frame association. In order to 
define a spinor structure our newly defined frame association should 
be consistent with the definition~\mythedefinition{2.6}. Let's prove
its consistence. For this purpose we consider two positively polarized
right orthonormal frames $\boldsymbol\Upsilon_0,\,\boldsymbol\Upsilon_1,
\,\boldsymbol\Upsilon_2,\,\boldsymbol\Upsilon_3$ and $\tilde{\boldsymbol
\Upsilon}_0,\,\tilde{\boldsymbol\Upsilon}_1,\,\tilde{\boldsymbol
\Upsilon}_2,\,\tilde{\boldsymbol\Upsilon}_3$ of the metric $\bold g$
whose domains are overlapping. They are related to each other with
the relationships \mythetag{2.2}, where $S$ and $T$ are mutually inverse
transition matrices. Applying the operator $\bold F$ to both sides of
the relationships \mythetag{2.2}, we get
$$
\xalignat 2
&\hskip -2em
\bold F(\tilde{\boldsymbol\Upsilon}_i)
=\sum^3_{j=0}S^j_i\,\bold F(\boldsymbol\Upsilon_j),
&&\bold F(\boldsymbol\Upsilon_i)
=\sum^3_{j=0}T^j_i\,\bold F(\tilde{\boldsymbol\Upsilon}_j).
\mytag{4.6}
\endxalignat
$$
The formulas \mythetag{4.6} mean that if we produce 
two positively polarized right orthonormal frames $\hat{\boldsymbol
\Upsilon}_0,\,\hat{\boldsymbol\Upsilon}_1,
\,\hat{\boldsymbol\Upsilon}_2,\,\hat{\boldsymbol\Upsilon}_3$ 
and $\hat{\Tilde{\boldsymbol\Upsilon}}_0,\,\hat{\Tilde{
\boldsymbol\Upsilon}}_1,\,\hat{\Tilde{\boldsymbol\Upsilon}}_2,
\,\hat{\Tilde{\boldsymbol\Upsilon}}_3$ of the metric $\hat\bold g$
by applying the deformation operator $\bold F$ to the frames 
$\boldsymbol\Upsilon_0,\,\boldsymbol\Upsilon_1,
\,\boldsymbol\Upsilon_2,\,\boldsymbol\Upsilon_3$ and $\tilde{\boldsymbol
\Upsilon}_0,\,\tilde{\boldsymbol\Upsilon}_1,\,\tilde{\boldsymbol
\Upsilon}_2,\,\tilde{\boldsymbol\Upsilon}_3$, then the deformed
frames will be related to each other by the same matrices
$S$ and $T$ in \mythetag{4.6} as the original frames are related
in \mythetag{2.2}:
$$
\xalignat 2
&\hskip -2em
\hat{\Tilde{\boldsymbol\Upsilon}}_i
=\sum^3_{j=0}S^j_i\,\hat{\boldsymbol\Upsilon}_j,
&&\hat{\boldsymbol\Upsilon}_i
=\sum^3_{j=0}T^j_i\,\hat{\Tilde{\boldsymbol\Upsilon}}_j.
\mytag{4.7}
\endxalignat
$$
Due to \mythetag{4.7} we can extend the diagram \mythetag{4.5} 
in the following way:
$$
\vcenter to 140pt{\vss\includegraphics{MetDef02.eps}\vss}
\mytag{4.8}
$$
Since the matrices $S$ and $T$ in the right hand side of the diagram
\mythetag{4.8} are the same as in its left hand side, the diagram
is commutative in whole, provided it is commutative in the absence
of the red arrows in it. The latter fact follows from the 
definition~\mythedefinition{2.6}. Moreover, the matrices 
$S$ and $T$ are related to the matrices $\goth S$ and $\goth T$
through the group homomorphism \mythetag{2.1}:
$$
\xalignat 2
&\hskip -2em
S=\phi(\goth S),
&&T=\phi(\goth T).
\endxalignat
$$
Therefore, the frame associations represented by the red horizontal
arrows on the diagram \mythetag{4.8} are consistent and define a new
spinor structure. Note that the spinor frames of this new spinor 
structure coincide with those of the initial one. This means that
the spinor structures produced by the metrics $\bold g$ and $\hat
\bold g$ share the same spinor bundle $SM$.
\mytheorem{4.1} The bundle of Weyl spinors $SM$ is preserved under
the metric deformations of the form \mythetag{3.1}.
\endproclaim
     Note that the frames $\boldsymbol\Psi_1,\,\boldsymbol\Psi_2$
and $\tilde{\boldsymbol\Psi}_1,\,\tilde{\boldsymbol\Psi}_2$ in the
diagram \mythetag{4.8} are canonically associated with positively 
polarized right orthonormal frames \pagebreak for both metrics 
$\bold g$ and $\hat\bold g$. Therefore, they are orthonormal spinor 
frames in the sense of the definition~\mythedefinition{2.7} for both 
spinor structures. This fact means that the spinor metrics $\bold d$ 
and $\hatboldd$ of both spinor structures are represented by the same
matrix \mythetag{2.11} in the same set of frames whose domains cover 
the space-time manifold. As a conclusion we get the following theorem.
\mytheorem{4.2} The skew-symmetric spinor metric $\bold d$ is
preserved under the metric deformations of the form \mythetag{3.1},
i\.\,e\. $\hatboldd=\bold d$.
\endproclaim
     Now let's consider the Infeld-van der Waerden field $\bold G$
for the non-deformed metric $\bold g$. It is represented by the
matrices \mythetag{2.14} in the frame pair 
$$
\hskip -2em
\boldsymbol\Upsilon_0,\,\boldsymbol\Upsilon_1,
\,\boldsymbol\Upsilon_2,\,\boldsymbol\Upsilon_3
\qquad\to\qquad\boldsymbol\Psi_1,\,\boldsymbol\Psi_2.
\mytag{4.9}
$$
Using the frame deformation formula \mythetag{4.1}, we can transform 
its components from the initial frame pair \mythetag{4.9} to the 
deformed frame pair
$$
\hskip -2em
\hat{\boldsymbol\Upsilon}_0,\,\hat{\boldsymbol\Upsilon}_1,
\,\hat{\boldsymbol\Upsilon}_2,\,\hat{\boldsymbol\Upsilon}_3
\qquad\to\qquad\boldsymbol\Psi_1,\,\boldsymbol\Psi_2.
\mytag{4.10}
$$
As a result for the components of $\bold G$ in the frame pair 
\mythetag{4.10} we get
$$
G^{i\kern 0.5pt\bar i}_k
=\sum^3_{q=0}F^q_k\,\sigma^{i\kern 0.5pt\bar i}_q.
\mytag{4.11}
$$
In the frame pair \mythetag{4.10} the components of the deformed 
Infeld-van der Waerden field $\hat\bold G$ are given by the Pauli 
matrices: $\hat G^{i\kern 0.5pt\bar i}_q=\sigma^{i\kern 0.5pt
\bar i}_q$. Applying this fact to the formula \mythetag{4.11}, we 
derive the following relationship:
$$
G^{i\kern 0.5pt\bar i}_k
=\sum^3_{q=0}F^q_k\,\hat G^{i\kern 0.5pt\bar i}_q.
\mytag{4.12}
$$
\mytheorem{4.3} The Infeld-van der Waerden field $\bold G$ is
transformed according to the formula \mythetag{4.12} under the 
metric deformations of the form \mythetag{3.1}.
\endproclaim
\noindent Now we apply the expansion \mythetag{3.2} to \mythetag{4.12}.
As a result, we get
$$
\hat G^{i\kern 0.5pt\bar i}_k
=G^{i\kern 0.5pt\bar i}_k
-\sum^3_{q=0}\varepsilon\,f^q_k
\,G^{i\kern 0.5pt\bar i}_q+\dots\,.
\mytag{4.13}
$$
Combining \mythetag{4.13} with the formula \mythetag{3.9} we derive the
following formula for the variation of the Infeld-van der Waerden field
$\bold G$:
$$
\hskip -2em
\delta G^{i\kern 0.5pt\bar i}_k
=-\frac{1}{2}\sum^3_{q=0}\sum^3_{p=0}g^{p\kern 0.3pt q}
\,G^{i\kern 0.5pt\bar i}_q\,\varepsilon\,h_{pk}.
\mytag{4.14}
$$
It is preferable to express $\delta G^{i\kern 0.5pt\bar i}_k$ through
$h^{ij}$ with upper indices. For this reason we transform the
formula \mythetag{4.14} as follows:
$$
\pagebreak 
\hskip -2em
\delta G^{i\kern 0.5pt\bar i}_k
=-\frac{1}{2}\sum^3_{q=0}\sum^3_{p=0}g_{pk}
\,G^{i\kern 0.5pt\bar i}_q\,\varepsilon\,h^{p\kern 0.3pt q}.
\mytag{4.15}
$$
Apart from \mythetag{4.15}, due to the theorem~\mythetheorem{4.2} we 
have the formula
$$
\hskip -2em
\delta d_{ij}=0.
\mytag{4.16}
$$
The formulas \mythetag{4.15} and \mythetag{4.16} describe completely 
the variations of the basic spin-tensorial fields \mythetag{2.15} of 
the bundle of Weyl spinors under the deformations of metric given by 
the formula \mythetag{3.1}. According the theorem~\mythetheorem{4.1},
the bundle $SM$ itself is invariant under these metric deformations.
\head
5. The bundle of Dirac spinors.
\endhead
     The bundle of Dirac spinors $DM$ over a space-time manifold $M$
is constructed as the direct sum of the bundle of Weyl spinors and its
Hermitian conjugate bundle:
$$
\hskip -2em
DM=SM\oplus S^{\sssize\dagger}\!M
\mytag{5.1}
$$
(see \mycite{3} for more details). Applying the 
theorem~\mythetheorem{4.1} to the expansion \mythetag{5.1} we
immediately derive the following theorem.
\mytheorem{5.1} The bundle of Dirac spinors $DM$ is preserved under
the metric deformations of the form \mythetag{3.1}.
\endproclaim
Like in the case of Weyl spinors, the structure of the Dirac bundle
$DM$ is described in terms of associated frame pairs and in terms of
basic spin-tensorial fields. There are four types of associated frame 
pairs 
$$
\hskip -2em
\aligned
&\boxit{Canonically orthonormal}{chiral frames}\leftarrow
\boxit{Positively polarized}{right orthonormal frames}\\
&\boxit{$P$-reverse}{antichiral frames}\leftarrow
\boxit{Positively polarized}{left orthonormal frames}\\
&\boxit{$T$-reverse}{antichiral frames}\leftarrow
\boxit{Negatively polarized}{right orthonormal frames}\\
&\boxit{$PT$-reverse}{chiral frames}\leftarrow
\boxit{Negatively polarized}{left orthonormal frames}
\endaligned
\mytag{5.2}
$$
and there are four basic spin-tensorial fields
$$
\vcenter{\hsize 10cm
\offinterlineskip\settabs\+\indent
\vtrule
\hskip 1.2cm &\vtrule % Quantity
\hskip 5.2cm &\vtrule % Unit
\hskip 2.8cm &\vtrule % Relation
\cr\hrule 
\+\vtrule
\hfill\,Symbol\hfill&\vtrule
\hfill Name\hfill &\vtrule
\hfill Spin-tensorial\hfill &\vtrule\cr
\vskip -0.2cm
\+\vtrule
\hfill &\vtrule
\hfill \hfill&\vtrule
\hfill type\hfill&\vtrule\cr\hrule
\+\vtrule
\hfill $\bold d$\hfill&\vtrule
\hfill Skew-symmetric metric tensor\hfill&\vtrule
\hfill $(0,2|0,0|0,0)$\hfill&\vtrule\cr\hrule
\+\vtrule
\hfill$\bold H$\hfill&\vtrule
\hfill Chirality operator\hfill&\vtrule
\hfill $(1,1|0,0|0,0)$\hfill&\vtrule\cr\hrule
\+\vtrule
\hfill$\bold D$\hfill&\vtrule
\hfill Dirac form\hfill&\vtrule
\hfill $(0,1|0,1|0,0)$\hfill&\vtrule\cr\hrule
\+\vtrule
\hfill$\boldsymbol\gamma$\hfill&\vtrule
\hfill Dirac $\gamma$-field\hfill&\vtrule
\hfill $(1,1|0,0|0,1)$\hfill&\vtrule\cr\hrule
}\quad
\mytag{5.3}
$$
in the bundle of Dirac spinors $DM$. Associated frame pairs
of the first type in \mythetag{5.2} are produced directly from
associated frame pairs of the bundle of Weyl spinors in
\mythetag{2.13}. They are sufficient for our purposes in this
paper. Associated frame pairs of the three other types are 
produced from associated frame pairs of the first
type by means of the $P$ and $T$-reversion procedures. We do
not consider them here referring the reader to the paper 
\mycite{3}.\par
     Let's consider some associated frame pair of the first type
for the bundle of Dirac spinors $DM$. It is composed of two frames:
$$
\hskip -2em
\boldsymbol\Upsilon_0,\,\boldsymbol\Upsilon_1,
\,\boldsymbol\Upsilon_2,\,\boldsymbol\Upsilon_3
\qquad\to\qquad\boldsymbol\Psi_1,\,\boldsymbol\Psi_2,
\,\boldsymbol\Psi_3,\,\boldsymbol\Psi_4.
\mytag{5.4}
$$
Here $\boldsymbol\Upsilon_0,\,\boldsymbol\Upsilon_1,
\,\boldsymbol\Upsilon_2,\,\boldsymbol\Upsilon_3$ is a positively
polarized right orthonormal frame of the tangent bundle $TM$. The 
second frame $\boldsymbol\Psi_1,\,\boldsymbol\Psi_2,\,\boldsymbol
\Psi_3,\,\boldsymbol\Psi_4$ in \mythetag{5.4} is a canonically 
orthonormal chiral frame of $DM$. Its first two spinor fields
$\boldsymbol\Psi_1$ and $\boldsymbol\Psi_2$ form an orthonormal
frame of the bundle of Weyl spinors $SM$ (see \mythetag{4.9} and
the definition~\mythedefinition{2.7} above). The second two fields
$\boldsymbol\Psi_3$ and $\boldsymbol\Psi_4$ belong to
$S^{\sssize\dagger}\!M$ in \mythetag{5.1}. They are semilinear
functional dual to $\boldsymbol\Psi_1$ and $\boldsymbol\Psi_2$
in the sense of the following relationships:
$$
\xalignat 2
&\hskip -2em
\boldsymbol\Psi_3(\boldsymbol\Psi_1)=1,
&&\boldsymbol\Psi_3(\boldsymbol\Psi_2)=0,\\
\vspace{-1ex}
\mytag{5.5}\\
\vspace{-1ex}
&\hskip -2em
\boldsymbol\Psi_4(\boldsymbol\Psi_1)=0,
&&\boldsymbol\Psi_4(\boldsymbol\Psi_2)=1.
\endxalignat
$$
The relationships \mythetag{5.5} fix $\boldsymbol\Psi_3$ and
$\boldsymbol\Psi_4$ uniquely, provided $\boldsymbol\Psi_1$ and 
$\boldsymbol\Psi_2$ are fixed. For this reason we can transform
the diagram \mythetag{4.5} as follows:
$$
\vcenter to 56pt{\vss\includegraphics{MetDef03.eps}\vss}
\mytag{5.6}
$$\par
    Once some associated frame pair \mythetag{5.4} of the first
type is fixed, the basic spin tensorial fields $\bold d$, 
$\bold H$, $\bold D$, and $\boldsymbol\gamma$ \mythetag{5.3}
are introduced by the following definitions.
\mydefinition{5.1} The skew-symmetric metric tensor $\bold d$ is
a spin-tensorial field of the type $(0,2|0,0|0,0)$ given by the
matrix 
$$
\hskip -2em
d_{ij}=\Vmatrix 0 & 1 & 0 & 0\\-1 & 0 & 0 & 0\\
0 & 0 & 0 & -1\\0 & 0 & 1 & 0\endVmatrix
\mytag{5.7}
$$
in any canonically orthonormal chiral frame of the Dirac bundle $DM$.
\enddefinition
\mydefinition{5.2} The chirality operator $\bold H$ is a spin-tensorial
field of the type $(1,1|0,0|0,0)$ given by the matrix 
$$
H^i_j=\Vmatrix 1 & 0 & 0 & 0\\0 & 1 & 0 & 0\\
0 & 0 & -1 & 0\\0 & 0 & 0 & -1\endVmatrix
\mytag{5.8}
$$
in any canonically orthonormal chiral frame of the Dirac bundle $DM$.
\enddefinition
\mydefinition{5.3} The Dirac form $\bold D$ is a spin-tensorial
field of the type $(0,1|0,1|0,0)$ given by the matrix 
$$
\hskip -2em
D_{i\bar j}=\Vmatrix 0 & 0 & 1 & 0\\0 & 0 & 0 & 1\\
1 & 0 & 0 & 0\\0 & 1 & 0 & 0\endVmatrix
\mytag{5.9}
$$
in any canonically orthonormal chiral frame of the Dirac bundle $DM$.
\enddefinition
\mydefinition{5.4} The Dirac $\gamma$-filed is a spin-tensorial 
field of the type $(1,1|0,0|0,1)$ given by the Dirac matrices  
$$
\xalignat 2
&\hskip -2em
\gamma^a_{b0}=m^a_{b0}=\Vmatrix 0&0&1&0\\0&0&0&1\\1&0&0&0
\\0&1&0&0\endVmatrix,
&&\gamma^a_{b1}=m^a_{b1}=\Vmatrix 0&0&0&1\\0&0&1&0\\0&-1&0&0
\\-1&0&0&0\endVmatrix,
\quad\\
\vspace{-1.5ex}
&&&\mytag{5.10}\\
\vspace{-1.5ex}
&\hskip -2em
\gamma^a_{b2}=m^a_{b2}=\Vmatrix 0&0&0&-i\\0&0&i&0\\0&i&0&0
\\-i&0&0&0\endVmatrix,
&&\gamma^a_{b3}=m^a_{b3}=\Vmatrix 0&0&1&0\\0&0&0&-1\\-1&0&0&0
\\0&1&0&0\endVmatrix
\quad
\endxalignat
$$
in any frame pair composed by a positively polarized right orthonormal
frame in $TM$ and its associated canonically orthonormal chiral frame 
in $DM$.
\enddefinition
\head
6. Deformation of the basic spin-tensorial fields.
\endhead
     Note that two positively polarized 
right orthonormal frames $\boldsymbol\Upsilon_0,\,\boldsymbol
\Upsilon_1,\,\boldsymbol\Upsilon_2,\,\boldsymbol\Upsilon_3$ and
$\hat{\boldsymbol\Upsilon}_0,\,\hat{\boldsymbol\Upsilon}_1,
\,\hat{\boldsymbol\Upsilon}_2,\,\hat{\boldsymbol\Upsilon}_3$
for the initial and deformed metrics on the diagram \mythetag{5.6} 
share the same canonically orthonormal chiral frame $\boldsymbol
\Psi_1,\,\boldsymbol\Psi_2,\,\boldsymbol\Psi_3,\,\boldsymbol\Psi_4$.
Note also that the spin-tensorial fields $\bold d$, $\bold H$, and
$\bold D$ are defined through spinor frames only (see 
definitions~\mythedefinition{5.1}, \mythedefinition{5.2}, and 
\mythedefinition{5.3} above). Therefore, the components of
$\bold d$, $\bold H$, and $\bold D$ coincide with the components
of $\hatboldd$, $\hat\bold H$, and $\hatboldD$ in the spinor 
frame $\boldsymbol\Psi_1,\,\boldsymbol\Psi_2,\,
\boldsymbol\Psi_3,\,\boldsymbol\Psi_4$ shared by two frame pairs
on the diagram \mythetag{5.6}. They are given by the matrices 
\mythetag{5.7}, \mythetag{5.8}, and \mythetag{5.9} respectively.
This result is expressed by the equalities
$$
\xalignat 3
&\bold d=\hatboldd,
&&\bold H=\hat\bold H,
&&\bold D=\hatboldD
\qquad
\mytag{6.1}
\endxalignat
$$ 
and formulated verbally in the following theorem.
\mytheorem{6.1} The skew-symmetric spinor metric\/ $\bold d$,
the chirality operator\/ $\bold H$ and the Dirac form\/ $\bold D$ 
are preserved under the metric deformations of the form 
\mythetag{3.1}.
\endproclaim
\noindent The formulas \mythetag{6.1} and the theorem~\mythetheorem{6.1} 
can be expressed in terms of variations:
$$
\xalignat 3
&\delta d_{ij}=0,
&&\delta H^i_j=0,
&&\delta D_{i\bar j}=0.
\qquad
\mytag{6.2}
\endxalignat
$$\par 
Now let's proceed to the Dirac $\gamma$-field. In this case $\boldsymbol
\gamma$ and $\hatboldsymbolgamma$ are different. The $\gamma$-field for
the non-deformed metric $\bold g$ is given by the $m$-matrices \mythetag{5.10}
in the non-deformed frame pair $\boldsymbol\Upsilon_0,
\,\boldsymbol\Upsilon_1,\,\boldsymbol\Upsilon_2,\,\boldsymbol\Upsilon_3
\to\boldsymbol\Psi_1,\,\boldsymbol\Psi_2,\,\boldsymbol\Psi_3,
\,\boldsymbol\Psi_4$ shown by the black right arrow on the diagram 
\mythetag{5.6}. Using the transition formulas \mythetag{4.1}, we can 
calculate its components in the deformed frame pair shown by the red 
arrow:
$$
\hskip -2em
\gamma^a_{bk}=\sum^3_{q=0}F^q_k\,m^a_{b\kern 0.5pt q}.
\mytag{6.3}
$$ 
This formula is analogous to \mythetag{4.11}. Now if we remember that 
the deformed $\gamma$-field is given by the $m$-matrices \mythetag{5.10}
in the deformed frame pair, we can transform \mythetag{6.3} to the 
formula analogous to \mythetag{4.12}:
$$
\hskip -2em
\gamma^a_{bk}=\sum^3_{q=0}F^q_k\,\hatgamma^a_{b\kern 0.5pt q}.
\mytag{6.4}
$$ 
\mytheorem{6.2} The Dirac $\gamma$-field is transformed according 
to the formula \mythetag{6.4} under the metric deformations of the 
form \mythetag{3.1}.
\endproclaim
\noindent Let's apply the expansion \mythetag{3.2} to $F^q_k$ in 
\mythetag{6.4}. As a result we obtain
$$
\hskip -2em
\hatgamma^a_{bk}=\gamma^a_{bk}-\sum^3_{q=0}\varepsilon\,f^q_k
\,\gamma^a_{b\kern 0.5pt q}+\dots\,.
\mytag{6.5}
$$
The formula \mythetag{6.5} is analogous to \mythetag{4.13}. Applying
\mythetag{3.9} to \mythetag{6.5}, we get
$$
\hskip -2em
\delta\gamma^a_{bk}=-\frac{1}{2}\sum^3_{q=0}\sum^3_{p=0}g^{p\kern 0.3pt q}
\,\gamma^a_{b\kern 0.5pt q}\,\varepsilon\,h_{pk}.
\mytag{6.6}
$$
Raising indices of $h_{pk}$ in \mythetag{6.6}, we can transform this formula
as follows:
$$
\hskip -2em
\delta\gamma^a_{bk}=-\frac{1}{2}\sum^3_{q=0}\sum^3_{p=0}g_{pk}
\,\gamma^a_{b\kern 0.5pt q}\,\varepsilon\,h^{p\kern 0.3pt q}.
\mytag{6.7}
$$
The formulas \mythetag{6.6} and \mythetag{6.7} are analogous to the formulas
\mythetag{4.14} and \mythetag{4.15} respectively. The formula \mythetag{6.7}
is complementary to \mythetag{6.2}.
\head
7. Deformation of the metric connection.
\endhead
     Each metric $\bold g$ of the space-time manifold $M$ produces the 
torsion-free metric connection $\Gamma$ in $TM$. This connection has 
extensions to the spinor bundles $SM$ and $DM$. In this paper we consider 
the extension $(\Gamma,\Alpha,\bar{\Alpha})$ of the metric connection 
$\Gamma$ to the bundle of Dirac spinors $DM$ only because it is preferably 
used in particle physics. Once some frame pair $\boldsymbol\Upsilon_0,
\,\boldsymbol\Upsilon_1,\,\boldsymbol\Upsilon_2,\,\boldsymbol\Upsilon_3$ 
and $\boldsymbol\Psi_1,\,\boldsymbol\Psi_2,\,\boldsymbol\Psi_3,
\,\boldsymbol\Psi_4$ is chosen, the components of the metric connection 
$(\Gamma,\Alpha,\bar{\Alpha})$ are given by explicit formulas. For its 
$\Gamma$-components we use the formula
$$
\hskip -2em
\gathered
\Gamma^k_{ij}=\sum^3_{r=0}\frac{g^{\kern 0.5pt kr}}{2}
\left(L_{\boldsymbol\Upsilon_i}\!(g_{rj})
+L_{\boldsymbol\Upsilon_j}\!(g_{i\kern 0.5pt r})
-L_{\boldsymbol\Upsilon_r}\!(g_{ij})\right)+\\
+\,\frac{c^{\,k}_{ij}}{2}
-\sum^3_{r=0}\sum^3_{s=0}\frac{c^{\,s}_{i\kern 0.5pt r}}{2}\,g^{kr}
\,g_{sj}-\sum^3_{r=0}\sum^3_{s=0}\frac{c^{\,s}_{j\kern 0.5ptr}}{2}
\,g^{kr}\,g_{s\kern 0.5pt i}
\endgathered
\mytag{7.1}
$$
taken from \mycite{4}. Here $L_{\boldsymbol\Upsilon_i}$, $L_{\boldsymbol
\Upsilon_j}$, and $L_{\boldsymbol\Upsilon_r}$ are the derivatives along 
the vectors $\boldsymbol\Upsilon_i$, $\boldsymbol\Upsilon_j$, and 
$\boldsymbol\Upsilon_r$ respectively, while $c^{\,k}_{ij}$ are the 
commutation coefficients taken from the commutation relationships 
\mythetag{2.3}.\par
     According to the theorem~\mythetheorem{4.1}, the bundle of Dirac 
spinors $DM$ is preserved under the deformations of metric. For this reason 
we fix a pair of frames $\boldsymbol\Upsilon_0,\,\boldsymbol\Upsilon_1,
\,\boldsymbol\Upsilon_2,\,\boldsymbol\Upsilon_3$ and $\boldsymbol\Psi_1,
\,\boldsymbol\Psi_2,\,\boldsymbol\Psi_3,\,\boldsymbol\Psi_4$ and assume 
these frames to be unchanged under the deformation of metric. Moreover,
we assume these frames to be canonically associated to each other 
according to the first line in the diagram \mythetag{5.2}. This means
that $\boldsymbol\Upsilon_0,\,\boldsymbol\Upsilon_1,\,\boldsymbol\Upsilon_2,
\,\boldsymbol\Upsilon_3$ is a positively polarized right orthonormal frame
with respect to the initial metric $\bold g$, while $\boldsymbol\Psi_1,
\,\boldsymbol\Psi_2,\,\boldsymbol\Psi_3,\,\boldsymbol\Psi_4$ is its
associated canonically orthonormal chiral frame in $DM$. Under these 
assumptions we have 
$$
\hskip -2em
\delta c^{\,k}_{ij}=0.
\mytag{7.2}
$$
Moreover, under these assumptions the components of the initial metric 
$\bold g$ are constants, i\.\,e\. $g_{ij}$ in \mythetag{3.7} are constants
taken from the matrix \mythetag{2.5}, while $\hat g_{ij}\neq\const$ and
$\delta g_{ij}$ in \mythetag{3.11} are also not constants. From 
$g_{ij}=\const$ we derive
$$
\hskip -2em
L_{\boldsymbol\Upsilon_k}\!(g_{ij})=0.
\mytag{7.3}
$$
Applying \mythetag{7.3} to the formula \mythetag{7.1}, we simplify it as 
follows:
$$
\hskip -2em
\gathered
\Gamma^k_{ij}=\frac{c^{\,k}_{ij}}{2}
-\sum^3_{r=0}\sum^3_{s=0}\frac{c^{\,s}_{i\kern 0.5pt r}}{2}\,g^{kr}
\,g_{sj}-\sum^3_{r=0}\sum^3_{s=0}\frac{c^{\,s}_{j\kern 0.5ptr}}{2}
\,g^{kr}\,g_{s\kern 0.5pt i}.
\endgathered
\mytag{7.4}
$$
Then we apply \mythetag{3.11}, \mythetag{3.18}, and \mythetag{7.2} to the
formula \mythetag{7.1}. As a result we get
$$
\hskip -2em
\gathered
\delta\Gamma^k_{ij}=-\sum^3_{r=0}\varepsilon\,\frac{g^{\kern 0.5pt kr}}{2}
\left(L_{\boldsymbol\Upsilon_i}\!(h_{rj})
+L_{\boldsymbol\Upsilon_j}\!(h_{i\kern 0.5pt r})
-L_{\boldsymbol\Upsilon_r}\!(h_{ij})\right)-\\
-\sum^3_{r=0}\sum^3_{s=0}\varepsilon\,\frac{c^{\,s}_{i\kern 0.5pt r}}{2}
\,h^{kr}\,g_{sj}-\sum^3_{r=0}\sum^3_{s=0}\varepsilon\,\frac{c^{\,s}_{j
\kern 0.5ptr}}{2}\,h^{kr}\,g_{s\kern 0.5pt i}\,+\\
+\sum^3_{r=0}\sum^3_{s=0}\varepsilon\,\frac{c^{\,s}_{i\kern 0.5pt r}}{2}
\,g^{kr}\,h_{sj}+\sum^3_{r=0}\sum^3_{s=0}\varepsilon\,\frac{c^{\,s}_{j
\kern 0.5ptr}}{2}\,g^{kr}\,h_{s\kern 0.5pt i}.
\endgathered
\mytag{7.5}
$$
The derivatives $L_{\boldsymbol\Upsilon_i}\!(h_{rj})$, $L_{\boldsymbol
\Upsilon_j}\!(h_{i\kern 0.5pt r})$, and $L_{\boldsymbol\Upsilon_r}
\!(h_{ij})$ in \mythetag{7.5} should be expressed through covariant 
derivatives. For this purpose we use the formula
$$
\hskip -2em
\nabla_{\!i}h_{jk}=L_{\boldsymbol\Upsilon_i}\!(h_{jk})-\sum^3_{s=0}
\Gamma^s_{ij}\,h_{sk}-\sum^3_{s=0}\Gamma^s_{ik}\,h_{js}.
\mytag{7.6}
$$
Applying \mythetag{7.6} to \mythetag{7.5} we obtain the following formula:
$$
\hskip -2em
\delta\Gamma^k_{ij}=-\sum^3_{r=0}\varepsilon\,\frac{g^{\kern 0.5pt kr}}{2}
\left(\nabla_{\!i}h_{rj}+\nabla_{\!j}h_{i\kern 0.5pt r}
-\nabla_{\!r}h_{ij}\right).
\mytag{7.7}
$$
This formula is well-known. Another proof of the formula \mythetag{7.7}
can be found in \mycite{1}.\par
     Now let's consider the $\Alpha$-components of the metric connection
$(\Gamma,\Alpha,\bar{\Alpha})$. They are given by the following formula 
taken from \mycite{5}:
$$
\gathered
\Alpha^a_{ib}
=\sum^4_{\alpha=1}\sum^4_{\beta=1}\frac{
L_{\boldsymbol\Upsilon_{\!i}}(d_{\alpha\beta})
\,d^{\kern 0.5pt \beta\kern 0.5pt\alpha}}{8}\,\delta^a_b
-\sum^4_{\alpha=1}\sum^4_{\beta=1}\sum^4_{d=1}\frac{
L_{\boldsymbol\Upsilon_{\!i}}(d_{\alpha\beta})
\,d^{\kern 0.5pt \beta\kern 0.5pt d}\,H^\alpha_d}{8}\,H^a_b\,-\\
-\sum^4_{c=1}\sum^4_{d=1}
\sum^4_{r=1}\frac{d_{\kern 0.5pt bc}
\,L_{\boldsymbol\Upsilon_i}(H^c_d)
\,H^d_r\,d^{\kern 0.5pt ra}}{4}\,+
\sum^3_{m=0}\sum^3_{n=0}\sum^4_{\alpha=1}
\frac{L_{\boldsymbol\Upsilon_i}(\gamma^{\,\alpha}_{b\kern 0.5pt m}
\,g^{mn})}{4}\,\times\\
\times\,\gamma^{\,a}_{\alpha\kern 0.5pt n}
+\sum^3_{m=0}\sum^3_{n=0}\sum^4_{\alpha=1}\sum^3_{s=0}
\frac{\gamma^{\,\alpha}_{b\kern 0.5pt m}\,\Gamma^n_{is}\,g^{ms}
\,\gamma^{\,a}_{\alpha\kern 0.5pt n}}{4}.
\endgathered\qquad
\mytag{7.8}
$$
Due to our special choice of frames $\boldsymbol\Upsilon_0,\,\boldsymbol
\Upsilon_1,\,\boldsymbol\Upsilon_2,\,\boldsymbol\Upsilon_3$ and $\boldsymbol
\Psi_1,\,\boldsymbol\Psi_2,\,\boldsymbol\Psi_3,\,\boldsymbol\Psi_4$ the
quantities $d_{\alpha\beta}$, $H^c_d$, $\gamma^{\,\alpha}_{b\kern 0.5pt m}$,
and $g^{mn}$ are constants. Their values are given by the formulas 
\mythetag{5.7}, \mythetag{5.8}, \mythetag{5.10}, and \mythetag{2.5}. For
this reason the formula \mythetag{7.8} reduces to
$$
\gathered
\Alpha^a_{ib}
=\sum^3_{m=0}\sum^3_{n=0}\sum^4_{\alpha=1}\sum^3_{s=0}
\frac{\gamma^{\,\alpha}_{b\kern 0.5pt m}\,\Gamma^n_{is}\,g^{ms}
\,\gamma^{\,a}_{\alpha\kern 0.5pt n}}{4}.
\endgathered\qquad
\mytag{7.9}
$$
The quantities $\Gamma^n_{is}$ in \mythetag{7.9} are given by 
\mythetag{7.4}. The formula \mythetag{7.9} is applicable 
only for the case of the initial non-deformed metric $\bold g$. 
Passing to its deformations \mythetag{3.1}, we should use the formula 
\mythetag{7.8} again. Taking into account \mythetag{6.2}, we get 
$$
\gathered
\delta\!\Alpha^a_{ib}
=\sum^3_{m=0}\sum^3_{n=0}\sum^4_{\alpha=1}
\frac{L_{\boldsymbol\Upsilon_i}(\delta\gamma^{\,\alpha}_{b\kern 0.5pt m}
\,g^{mn}+\gamma^{\,\alpha}_{b\kern 0.5pt m}\,\delta g^{mn})}{4}
\,\gamma^{\,a}_{\alpha\kern 0.5pt n}\,+\\
\vspace{1ex}
+\sum^3_{m=0}\sum^3_{n=0}\sum^4_{\alpha=1}\sum^3_{s=0}
\frac{\Gamma^n_{is}\,(\delta\gamma^{\,\alpha}_{b\kern 0.5pt m}\,g^{ms}+
\gamma^{\,\alpha}_{b\kern 0.5pt m}\,\delta g^{ms})
\,\gamma^{\,a}_{\alpha\kern 0.5pt n}}{4}\,+\\
\vspace{1ex}
+\sum^3_{m=0}\sum^3_{n=0}\sum^4_{\alpha=1}\sum^3_{s=0}
\frac{\gamma^{\,\alpha}_{b\kern 0.5pt m}\,\delta\Gamma^n_{is}\,g^{ms}
\,\gamma^{\,a}_{\alpha\kern 0.5pt n}+
\gamma^{\,\alpha}_{b\kern 0.5pt m}\,\Gamma^n_{is}\,g^{ms}
\,\delta\gamma^{\,a}_{\alpha\kern 0.5pt n}}{4}.
\endgathered
\quad
\mytag{7.10}
$$
Before transforming the formula \mythetag{7.10} in whole we consider some 
smaller subexpression in the right hand side of this formula. Applying 
the formulas \mythetag{3.18} and \mythetag{6.7} to this subexpression, 
we get
$$
\gathered
\sum^3_{m=0}(\delta\gamma^{\,\alpha}_{b\kern 0.5pt m}\,g^{mn}
+\gamma^{\,\alpha}_{b\kern 0.5pt m}\,\delta g^{mn})=
-\frac{1}{2}\sum^3_{q=0}\sum^3_{p=0}\sum^3_{m=0}g_{pm}
\,\gamma^{\,\alpha}_{b\kern 0.5pt q}\,\varepsilon\,h^{p\kern 0.3pt q}
\,g^{mn}\,+\\
+\sum^3_{m=0}\gamma^{\,\alpha}_{b\kern 0.5pt m}\,\varepsilon\,h^{mn}
=\frac{1}{2}\sum^3_{m=0}\gamma^{\,\alpha}_{b\kern 0.5pt m}\,\varepsilon
\,h^{mn}.
\endgathered
\quad
\mytag{7.11}
$$
The subexpression \mythetag{7.11} is differentiated in \mythetag{7.10}. 
Before substituting \mythetag{7.11} back into \mythetag{7.10} we express
the Lie derivative $L_{\boldsymbol\Upsilon_i}$ of it through its covariant
derivative:
$$
\hskip -2em
\gathered
\sum^3_{m=0}L_{\boldsymbol\Upsilon_i}(\delta\gamma^{\,\alpha}_{b
\kern 0.5pt m}\,g^{mn}+\gamma^{\,\alpha}_{b\kern 0.5pt m}\,\delta 
g^{mn})=\frac{\varepsilon}{2}\sum^3_{m=0}L_{\boldsymbol\Upsilon_i}
(\gamma^{\,\alpha}_{b\kern 0.5pt m}\,h^{mn})=\\
=\frac{1}{2}\sum^3_{m=0}\nabla_{\!i}(\gamma^{\,\alpha}_{b\kern 0.5pt m}
\,\varepsilon\,h^{mn})-\frac{1}{2}\sum^3_{m=0}\sum^4_{\theta=1}
A^\alpha_{i\kern 0.5pt\theta}\,\gamma^{\,\theta}_{b\kern 0.5pt m}\,
\varepsilon\,h^{mn}\,+\\
+\,\frac{1}{2}\sum^3_{m=0}\sum^4_{\theta=1}A^\theta_{ib}
\,\gamma^{\,\alpha}_{\theta\kern 0.5pt m}\,\varepsilon\,h^{mn}
-\,\frac{1}{2}\sum^3_{m=0}\sum^3_{s=0}\Gamma^n_{is}
\,\gamma^{\,\alpha}_{b\kern 0.5pt m}\,\varepsilon\,h^{ms}.
\endgathered
\mytag{7.12}
$$
Now, substituting \mythetag{7.12} into \mythetag{7.10}, we recall the 
following identities:
$$
\xalignat 5
&\nabla\bold d=0,
&&\nabla\bold H=0,
&&\nabla\bold D=0,
&&\nabla\boldsymbol\gamma=0,
&&\nabla\bold g=0.
\qquad\quad
\mytag{7.13}
\endxalignat
$$
The identities \mythetag{7.13} are known as the concordance conditions for 
the metric $\bold g$ and its metric connection $(\Gamma,\Alpha,\bar{\Alpha})$. 
The spin-tensorial fields $\bold d$, $\bold H$, $\bold D$, and $\boldsymbol
\gamma$ are treated as attributes of the metric $\bold g$. Thus, applying 
the identities \mythetag{7.13} to \mythetag{7.12} and substituting 
\mythetag{7.12} back into the formula \mythetag{7.10}, we get
$$
\gathered
\delta\!\Alpha^a_{ib}
=\sum^3_{m=0}\sum^3_{n=0}\sum^4_{\alpha=1}\left(
\frac{
\gamma^{\,\alpha}_{b\kern 0.5pt m}
\,\varepsilon\,\nabla_{\!i}h^{mn}
}{8}-\sum^4_{\theta=1}\frac{A^\alpha_{i\kern 0.5pt\theta}
\,\gamma^{\,\theta}_{b\kern 0.5pt m}\,\varepsilon\,h^{mn}}{8}\,
+\right.\\
\left.+\sum^4_{\theta=1}\,\frac{A^\theta_{i\kern 0.2pt b}
\,\gamma^{\,\alpha}_{\theta\kern 0.5pt m}\,\varepsilon\,h^{mn}}{8}
\right)\gamma^{\,a}_{\alpha\kern 0.5pt n}\,+\\
\vspace{1ex}
+\sum^3_{m=0}\sum^3_{n=0}\sum^4_{\alpha=1}\sum^3_{s=0}
\frac{\gamma^{\,\alpha}_{b\kern 0.5pt m}\,\delta\Gamma^n_{is}\,g^{ms}
\,\gamma^{\,a}_{\alpha\kern 0.5pt n}+
\gamma^{\,\alpha}_{b\kern 0.5pt m}\,\Gamma^n_{is}\,g^{ms}
\,\delta\gamma^{\,a}_{\alpha\kern 0.5pt n}}{4}.
\endgathered
\quad
\mytag{7.14}
$$
In the next step we apply the formula \mythetag{7.7} for $\delta\Gamma^n_{is}$
and the formula \mythetag{6.7} for $\delta\gamma^{\,a}_{\alpha\kern 0.5pt n}$
in order to transform the formula \mythetag{7.14}. As a result we get
$$
\gathered
\delta\!\Alpha^a_{ib}
=\sum^3_{m=0}\sum^3_{n=0}\sum^4_{\alpha=1}\sum^4_{\theta=1}\left(
\frac{A^\theta_{i\kern 0.2pt b}
\,\gamma^{\,\alpha}_{\theta\kern 0.5pt m}}{8}
-\frac{A^\alpha_{i\kern 0.5pt\theta}
\,\gamma^{\,\theta}_{b\kern 0.5pt m}}{8}\,
\right)\gamma^{\,a}_{\alpha\kern 0.5pt n}\,\varepsilon\,h^{mn}\,+\\
\vspace{1ex}
+\sum^3_{m=0}\sum^3_{n=0}\sum^4_{\alpha=1}\sum^3_{r=0}\sum^3_{s=0}
\frac{g_{in}\,\gamma^{\,\alpha}_{b\kern 0.5pt m}\,g^{rs}
\,\gamma^{\,a}_{\alpha\kern 0.5pt s}\,\varepsilon\,\nabla_{\!r}h^{mn}
}{8}\,-\\
\vspace{1ex}
-\sum^3_{m=0}\sum^3_{n=0}\sum^4_{\alpha=1}\sum^3_{r=0}\sum^3_{s=0}
\frac{g_{in}\,\gamma^{\,\alpha}_{b\kern 0.5pt s}\,g^{rs}
\,\gamma^{\,a}_{\alpha\kern 0.5pt m}\,\varepsilon\,\nabla_{\!r}h^{mn}
}{8}\,-\\
\vspace{1ex}
-\sum^3_{n=0}\sum^3_{m=0}\sum^3_{p=0}\sum^3_{q=0}\sum^4_{\alpha=1}
\sum^3_{s=0}
\frac{\gamma^{\,\alpha}_{b\kern 0.5pt p}\,\Gamma^q_{is}\,g^{p\kern 0.3pt s}
\,g_{mq}\,\gamma^a_{\alpha\kern 0.5pt n}
\,\varepsilon\,h^{m\kern 0.3pt n}}{8}.
\endgathered
\quad
\mytag{7.15}
$$
Now let's remember that $\nabla_{\!i}\gamma^{\,\alpha}_{b\kern 0.5pt m}=0$
and $\nabla_{\!i}\,g^{p\kern 0.5pt q}=0$. \pagebreak These identities are 
the coordinate forms of $\nabla\boldsymbol\gamma$ and $\nabla\bold g$ from 
\mythetag{7.13}. Expanding them, we get
$$
\align
&\nabla_{\!i}\gamma^{\,\alpha}_{b\kern 0.5pt m}=
L_{\boldsymbol\Upsilon_i}(\gamma^{\,\alpha}_{b\kern 0.5pt m})
+\sum^4_{\theta=1}A^\alpha_{i\kern 0.5pt\theta}
\,\gamma^{\,\theta}_{b\kern 0.5pt m}-
\sum^4_{\theta=1}A^\theta_{i\kern 0.2pt b}
\,\gamma^{\,\alpha}_{\theta\kern 0.5pt m}
-\sum^3_{s=0}\Gamma^s_{im}\,\gamma^{\,\alpha}_{\theta\kern 0.2pt s}=0,\\
&\nabla_{\!i}\,g^{p\kern 0.5pt q}=L_{\boldsymbol\Upsilon_i}(g^{p\kern 0.5pt
q})+\sum^3_{s=0}\Gamma^p_{is}\,g^{s\kern 0.1pt q}
+\sum^3_{s=0}\Gamma^q_{is}\,g^{p\kern 0.3pt s}=0.
\endalign
$$
Due to our special choice of frames $\boldsymbol\Upsilon_0,\,\boldsymbol
\Upsilon_1,\,\boldsymbol\Upsilon_2,\,\boldsymbol\Upsilon_3$ and $\boldsymbol
\Psi_1,\,\boldsymbol\Psi_2,\,\boldsymbol\Psi_3,\,\boldsymbol\Psi_4$ the
quantities $\gamma^{\,\alpha}_{b\kern 0.5pt m}$ and $g^{p\kern 0.5pt q}$
are constants. Therefore $L_{\boldsymbol\Upsilon_i}(\gamma^{\,\alpha}_{b
\kern 0.5pt m})=0$ and $L_{\boldsymbol\Upsilon_i}(g^{p\kern 0.5pt q})=0$.
As a result we obtain the following identities:
$$
\align
&\hskip -2em
\sum^4_{\theta=1}A^\theta_{i\kern 0.2pt b}
\,\gamma^{\,\alpha}_{\theta\kern 0.5pt m}
-\sum^4_{\theta=1}A^\alpha_{i\kern 0.5pt\theta}
\,\gamma^{\,\theta}_{b\kern 0.5pt m}
=-\sum^3_{s=0}\Gamma^s_{im}
\,\gamma^{\,\alpha}_{\theta\kern 0.2pt s},
\mytag{7.16}\\
&\hskip -2em
\sum^3_{s=0}\Gamma^q_{is}\,g^{p\kern 0.3pt s}
=-\sum^3_{s=0}\Gamma^p_{is}\,g^{s\kern 0.1pt q}.
\mytag{7.17}
\endalign
$$
Applying \mythetag{7.16} to the first term in \mythetag{7.15} and
applying \mythetag{7.17} to the last term in \mythetag{7.15}, we find 
that these terms cancel each other. This yields 
$$
\gathered
\delta\!\Alpha^a_{ib}
=\sum^3_{m=0}\sum^3_{n=0}\sum^4_{\alpha=1}\sum^3_{r=0}\sum^3_{s=0}
\frac{g_{in}\,\gamma^{\,\alpha}_{b\kern 0.5pt m}\,g^{rs}
\,\gamma^{\,a}_{\alpha\kern 0.5pt s}\,\varepsilon\,\nabla_{\!r}h^{mn}
}{8}\,-\\
\vspace{1ex}
-\sum^3_{m=0}\sum^3_{n=0}\sum^4_{\alpha=1}\sum^3_{r=0}\sum^3_{s=0}
\frac{g_{in}\,\gamma^{\,\alpha}_{b\kern 0.5pt s}\,g^{rs}
\,\gamma^{\,a}_{\alpha\kern 0.5pt m}\,\varepsilon\,\nabla_{\!r}h^{mn}
}{8}.
\endgathered
\quad
\mytag{7.18}
$$\par
     Passing to the barred $\Alpha$ components of the metric connection
$(\Gamma,\Alpha,\bar{\Alpha})$ remember that they are obtained as the 
complex conjugates of $\Alpha$ components:
$$
\hskip -2em
\bar{\Alpha}\vphantom{\Alpha}^{\bar a}_{i\kern 1pt \bar b}=
\overline{\Alpha^{\bar a}_{i\kern 1pt \bar b}}.
\mytag{7.19}
$$
From \mythetag{7.19} we immediately derive the formula
$$
\hskip -2em
\delta\bar{\Alpha}\vphantom{\Alpha}^{\bar a}_{i\kern 1pt \bar b}=
\overline{\delta\Alpha^{\bar a}_{i\kern 1pt \bar b}}.
\mytag{7.20}
$$
The formulas \mythetag{7.7}, \mythetag{7.18}, and \mythetag{7.20} describe 
completely the variation of the metric connection $(\Gamma,\Alpha,
\bar{\Alpha})$ under the metric deformations of the form \mythetag{3.1}.
\head
8. Energy-momentum tensor\\
of the massive spin 1/2 particle.
\endhead
     A single spin 1/2 particle with the mass $m$ in a space-time manifold 
$M$ is described by a spinor-valued $\boldsymbol\psi$-function satisfying 
the Dirac equation:
$$
\hskip -2em
i\,\hbar\,\sum^4_{b=1}\sum^3_{p=0}\sum^3_{q=0}
\gamma^{\kern 0.5pt a}_{b\kern 0.2pt p}
\,g^{p\kern 0.3pt q}\,\nabla_{\!q}
\psi^{\kern 0.5pt\lower 1.2pt\hbox{$\ssize b$}}
-m\,c\ \psi^{\kern 0.5pt\lower 1.2pt\hbox{$\ssize a$}}=0.
\mytag{8.1}
$$
The Dirac equation \mythetag{8.1} is derived from the following action
functional:
$$
\gathered
S_{\text{matter}}=i\,\hbar
\int\sum^4_{a=1}\sum^4_{\bar a=1}\sum^4_{b=1}\sum^3_{p=0}
\sum^3_{q=0}
D_{a\bar a}\,
\gamma^{\kern 0.5pt a}_{b\kern 0.2pt p}\,g^{p\kern 0.3pt q}
\,\frac{\overline{\psi^{\kern 0.5pt\bar a}}
\ \nabla_{\!q}\psi^{\kern 0.5pt\lower 1.2pt\hbox{$\ssize b$}}
-\psi^{\kern 0.5pt\lower 1.2pt\hbox{$\ssize b$}}\ \nabla_{\!q}
\overline{\psi^{\kern 0.5pt\bar a}}}{2}\,dV\,-\\
\vspace{2ex}
-\,m\,c\int\sum^4_{a=1}\sum^4_{\bar a=1}D_{a\bar a}\,
\overline{\psi^{\kern 0.5pt\bar a}}\,
\psi^{\kern 0.5pt\lower 1.2pt\hbox{$\ssize a$}}\,dV
\endgathered
\quad
\mytag{8.2}
$$
(see \mycite{6}). Our goal is to substitute the action functional 
\mythetag{8.2} into the formula \mythetag{1.4}. For this purpose 
we write the action \mythetag{8.2} formally as
$$
\hskip -2em
S_{\text{matter}}=\int L\,dV.
\mytag{8.3}
$$
The real scalar quantity $L$ in \mythetag{8.3} is called the
Lagrangian density. The variation of the action integral
\mythetag{8.3} is written as follows:
$$
\hskip -2em
\delta S_{\text{matter}}=\int\delta L\,dV
-\frac{1}{2}\int\sum^3_{i=0}\sum^3_{j=0}L\,g_{ij}
\,\varepsilon\,h^{ij}\,dV.
\mytag{8.4}
$$
The second integral in \mythetag{8.4} arises because the volume
element $dV$ depends on the metric $\bold g$ and changes under 
the metric deformations of the form \mythetag{3.1} (see more
details in \mycite{1}).\par
     Due to the formula \mythetag{8.4} we need to calculate $\delta L$. 
For this purpose we subdivide $L$ into subexpressions and calculate 
their variations separately. According to \mythetag{8.2}, the Lagrangian
density $L$ is a sum of two terms:
$$
\hskip -2em
L=L_{\text{kinetic}}+L_{\text{massive}}.
\mytag{8.5}
$$
The massive term in \mythetag{8.5} is given by the formula
$$
\hskip -2em
L_{\text{massive}}=-m\,c\sum^4_{a=1}\sum^4_{\bar a=1}D_{a\bar a}\,
\overline{\psi^{\kern 0.5pt\bar a}}\,
\psi^{\kern 0.5pt\lower 1.2pt\hbox{$\ssize a$}}.
\mytag{8.6}
$$
The variation of the massive term \mythetag{8.6} is equal to zero. Indeed,
we have 
$$
\xalignat 2
&\hskip -2em
\delta\psi^{\kern 0.5pt\lower 1.2pt\hbox{$\ssize a$}}=0,
&&\delta\overline{\psi^{\kern 0.5pt\bar a}}=0
\mytag{8.7}
\endxalignat
$$
since $\boldsymbol\psi$-function is treated as an independent parameter 
under the metric deformations \mythetag{3.1}. As for $D_{a\bar a}$ in
\mythetag{8.6}, for these parameters we have
$$
\hskip -2em
\delta D_{a\bar a}=0.
\mytag{8.8}
$$
The equality \mythetag{8.8} is derived from \mythetag{7.13}. The equalities
\mythetag{8.7} and \mythetag{8.8} lead to
$$
\hskip -2em
\delta L_{\text{massive}}=0.
\mytag{8.9}
$$\par
     Let's proceed to the kinetic term of the Lagrangian density
\mythetag{8.5}. According to \mythetag{8.2}, this term is given by 
the formula
$$
L_{\text{kinetic}}
=i\,\hbar\sum^4_{a=1}\sum^4_{\bar a=1}\sum^4_{b=1}\sum^3_{p=0}
\sum^3_{q=0}
D_{a\bar a}\,
\gamma^{\kern 0.5pt a}_{b\kern 0.2pt p}\,g^{p\kern 0.3pt q}
\,\frac{\overline{\psi^{\kern 0.5pt\bar a}}
\ \nabla_{\!q}\psi^{\kern 0.5pt\lower 1.2pt\hbox{$\ssize b$}}
-\psi^{\kern 0.5pt\lower 1.2pt\hbox{$\ssize b$}}\ \nabla_{\!q}
\overline{\psi^{\kern 0.5pt\bar a}}}{2}.
\quad
\mytag{8.10}
$$
Though due to \mythetag{8.7} $\psi^{\kern 0.5pt\lower 1.2pt\hbox{$\ssize 
b$}}$ and $\overline{\psi^{\kern 0.5pt\bar a}}$ in \mythetag{8.10}
are not sensitive to the metric deformations \mythetag{3.1}, their 
covariant derivatives are sensitive. Indeed, we have
$$
\hskip -2em
\aligned
&\nabla_{\!q}\psi^{\kern 0.5pt\lower 1.2pt\hbox{$\ssize b$}}
=L_{\boldsymbol\Upsilon_q}(\psi^{\kern 0.5pt\lower 1.2pt\hbox{$\ssize b$}})
+\sum^4_{\theta=1}\Alpha^b_{q\kern 0.4pt\theta}\,
\psi^{\kern 0.5pt\lower 1.2pt\hbox{$\ssize\theta$}},\\
&\nabla_{\!q}\overline{\psi^{\kern 0.5pt\bar a}}
=L_{\boldsymbol\Upsilon_q}(\overline{\psi^{\kern 0.5pt\bar a}})
+\sum^4_{\theta=1}\bar{\Alpha}\vphantom{\Alpha}^{\bar a}_{q\kern 0.4pt
\theta }\,\overline{\psi^{\kern 0.5pt\theta}}.
\endaligned
\mytag{8.11}
$$
From the equalities \mythetag{8.11} we immediately derive
$$
\xalignat 2
&\hskip -2em
\delta\nabla_{\!q}\psi^{\kern 0.5pt\lower 1.2pt\hbox{$\ssize b$}}
=\sum^4_{\theta=1}\delta\!\Alpha^b_{q\kern 0.4pt\theta}\,
\psi^{\kern 0.5pt\lower 1.2pt\hbox{$\ssize\theta$}},
&&\delta\nabla_{\!q}\overline{\psi^{\kern 0.5pt\bar a}}
=\sum^4_{\theta=1}\delta\bar{\Alpha}
\vphantom{\Alpha}^{\bar a}_{q\kern 0.4pt\theta }
\,\overline{\psi^{\kern 0.5pt\theta}}.
\quad
\mytag{8.12}
\endxalignat
$$
Apart from \mythetag{8.11}, there are two other terms 
$\gamma^{\kern 0.5pt a}_{b\kern 0.2pt p}$ and $g^{p\kern 0.3pt q}$
in the formula \mythetag{8.10} which are sensitive to the metric 
deformations \mythetag{3.1}. Their variations are given by the formulas 
\mythetag{6.7} and \mythetag{3.18}.\par
     Now let's calculate the variation of the kinetic term of the 
Lagrangian density \mythetag{8.5}. From the formula \mythetag{8.10},
applying \mythetag{8.12}, we derive
$$
\hskip -2em
\gathered
\delta L_{\text{kinetic}}
=i\,\hbar\sum^4_{a=1}\sum^4_{\bar a=1}\sum^4_{b=1}\sum^3_{p=0}
\sum^3_{q=0}D_{a\bar a}\,(
\delta\gamma^{\kern 0.5pt a}_{b\kern 0.2pt p}\,g^{p\kern 0.3pt q}
+\gamma^{\kern 0.5pt a}_{b\kern 0.2pt p}\,\delta g^{p\kern 0.3pt q})
\times\\
\vspace{2ex}
\times\,\frac{\overline{\psi^{\kern 0.5pt\bar a}}\,
\nabla_{\!q}\psi^{\kern 0.5pt\lower 1.2pt\hbox{$\ssize b$}}
-\psi^{\kern 0.5pt\lower 1.2pt\hbox{$\ssize b$}}\ \nabla_{\!q}
\overline{\psi^{\kern 0.5pt\bar a}}}{2}
+i\,\hbar\sum^4_{a=1}\sum^4_{\bar a=1}\sum^4_{b=1}\sum^3_{p=0}
\sum^3_{q=0}D_{a\bar a}\,\gamma^{\kern 0.5pt a}_{b\kern 0.2pt p}
\,g^{p\kern 0.3pt q}\,\times\\
\times\,\frac{1}{2}\sum^4_{\theta=1}
\left(\overline{\psi^{\kern 0.5pt\bar a}}\,
\delta\!\Alpha^b_{q\kern 0.4pt\theta}\,
\psi^{\kern 0.5pt\lower 1.2pt\hbox{$\ssize\theta$}}
-\psi^{\kern 0.5pt\lower 1.2pt\hbox{$\ssize b$}}
\, \delta\bar{\Alpha}
\vphantom{\Alpha}^{\bar a}_{q\kern 0.4pt\theta }
\,\overline{\psi^{\kern 0.5pt\theta}}
\right).
\endgathered
\quad
\mytag{8.13}
$$
Applying \mythetag{7.18} and \mythetag{7.20} to \mythetag{8.13}, we find 
that second term in \mythetag{8.13} is identically equal to zero. Therefore,
applying \mythetag{7.11} to \mythetag{8.13}, we obtain
$$
\gathered
\delta L_{\text{kinetic}}
=i\,\hbar\sum^4_{a=1}\sum^4_{\bar a=1}\sum^4_{b=1}\sum^3_{p=0}
\sum^3_{q=0}D_{a\bar a}\,\gamma^{\kern 0.5pt a}_{b\kern 0.2pt p}
\,\times\\
\vspace{2ex}
\kern 13em\times\,\frac{\overline{\psi^{\kern 0.5pt\bar a}}
\ \nabla_{\!q}\psi^{\kern 0.5pt\lower 1.2pt\hbox{$\ssize b$}}
-\psi^{\kern 0.5pt\lower 1.2pt\hbox{$\ssize b$}}
\ \nabla_{\!q}\overline{\psi^{\kern 0.5pt\bar a}}}{4}
\ \varepsilon\,h^{p\kern 0.3pt q}.
\endgathered
\qquad\quad
\mytag{8.14}
$$
Let's remember that $h^{p\kern 0.3pt q}$ is symmetric with respect to the
indices $p$ and $q$. Our next goal is to make other terms in \mythetag{8.14}
\pagebreak symmetric with respect to these indices. Applying the 
symmetrization procedure to \mythetag{8.14}, we remember \mythetag{8.9} which 
means that $\delta L=\delta L_{\text{kinetic}}$. Then the formula 
\mythetag{8.14} yields
$$
\gathered
\delta L
=i\,\hbar\sum^3_{p=0}\sum^3_{q=0}\!\left(\,\sum^4_{a=1}\sum^4_{\bar a=1}
\sum^4_{b=1}D_{a\bar a}\,
\frac{\gamma^{\kern 0.5pt a}_{b\kern 0.2pt p}
\ \overline{\psi^{\kern 0.5pt\bar a}}
\ \nabla_{\!q}\psi^{\kern 0.5pt\lower 1.2pt\hbox{$\ssize b$}}
+\gamma^{\kern 0.5pt a}_{b\kern 0.2pt q}
\ \overline{\psi^{\kern 0.5pt\bar a}}
\ \nabla_{\!p}\psi^{\kern 0.5pt\lower 1.2pt\hbox{$\ssize b$}}
}{8}\,-\right.\\
\vspace{2ex}
\qquad\qquad\left.-\sum^4_{a=1}\sum^4_{\bar a=1}\sum^4_{b=1}
D_{a\bar a}\,
\frac{\gamma^{\kern 0.5pt a}_{b\kern 0.2pt p}
\ \psi^{\kern 0.5pt\lower 1.2pt\hbox{$\ssize b$}}
\ \nabla_{\!q}\overline{\psi^{\kern 0.5pt\bar a}}
+\gamma^{\kern 0.5pt a}_{b\kern 0.2pt q}
\ \psi^{\kern 0.5pt\lower 1.2pt\hbox{$\ssize b$}}
\ \nabla_{\!p}\overline{\psi^{\kern 0.5pt\bar a}}
}{8}\,\right)\varepsilon\,h^{p\kern 0.3pt q}.
\endgathered
\qquad\quad
\mytag{8.15}
$$
Now, substituting \mythetag{8.15} into \mythetag{8.4} and applying 
the formula \mythetag{1.4}, we derive the explicit formula for the
components of the energy-momentum tensor
$$
\hskip -2em
\gathered
T_{ij}=i\hbar\,c\,\sum^4_{a=1}\sum^4_{\bar a=1}\sum^4_{b=1}
D_{a\bar a}\,
\frac{\gamma^{\kern 0.5pt a}_{b\kern 0.2pt i}
\ \overline{\psi^{\kern 0.5pt\bar a}}
\ \nabla_{\!j}\psi^{\kern 0.5pt\lower 1.2pt\hbox{$\ssize b$}}
+\gamma^{\kern 0.5pt a}_{b\kern 0.2pt j}
\ \overline{\psi^{\kern 0.5pt\bar a}}
\ \nabla_{\!i}\psi^{\kern 0.5pt\lower 1.2pt\hbox{$\ssize b$}}
}{4}\,-\\
-\,i\hbar\,c\,\sum^4_{a=1}\sum^4_{\bar a=1}\sum^4_{b=1}
D_{a\bar a}\,
\frac{\gamma^{\kern 0.5pt a}_{b\kern 0.2pt i}
\ \psi^{\kern 0.5pt\lower 1.2pt\hbox{$\ssize b$}}
\ \nabla_{\!j}\overline{\psi^{\kern 0.5pt\bar a}}
+\gamma^{\kern 0.5pt a}_{b\kern 0.2pt j}
\ \psi^{\kern 0.5pt\lower 1.2pt\hbox{$\ssize b$}}
\ \nabla_{\!i}\overline{\psi^{\kern 0.5pt\bar a}}
}{4}\,+\\
+\,i\hbar\,c\,\sum^4_{a=1}\sum^4_{\bar a=1}\sum^4_{b=1}
\sum^3_{p=0}\sum^3_{q=0}D_{a\bar a}
\,\frac{\gamma^{\kern 0.5pt a}_{b\kern 0.2pt p}\,g^{p\kern 0.3pt q}
\,\psi^{\kern 0.5pt\lower 1.2pt\hbox{$\ssize b$}}\ \nabla_{\!q}
\overline{\psi^{\kern 0.5pt\bar a}}}{2}\,g_{ij}\,-\\
-\,\,i\hbar\,c\,\sum^4_{a=1}\sum^4_{\bar a=1}\sum^4_{b=1}
\sum^3_{p=0}\sum^3_{q=0}D_{a\bar a}
\,\frac{\gamma^{\kern 0.5pt a}_{b\kern 0.2pt p}\,g^{p\kern 0.3pt q}
\overline{\psi^{\kern 0.5pt\bar a}}
\ \nabla_{\!q}\psi^{\kern 0.5pt\lower 1.2pt\hbox{$\ssize b$}}}{2}
\,g_{ij}\,+\\
\kern -80pt+\,m\,c^2\,\sum^4_{a=1}\sum^4_{\bar a=1}D_{a\bar a}\,
\overline{\psi^{\kern 0.5pt\bar a}}\,
\psi^{\kern 0.5pt\lower 1.2pt\hbox{$\ssize a$}}
\,g_{ij}.
\endgathered
\mytag{8.16}
$$
Though the $\boldsymbol\psi$-function is a complex-valued 
spinor function, the components of the energy-momentum tensor
\mythetag{8.16} are real. This fact is proved on the base of 
the following identity relating the components of the 
fields $\bold D$ and $\boldsymbol\gamma$:
$$
\hskip -2em
\sum^4_{a=1}\overline{D_{a\bar a}}\
\overline{\gamma^{\kern 0.5pt a}_{b\kern 0.2pt p}
\vphantom{D_{a\bar a}}}
=\sum^4_{a=1}D_{ab}\ \gamma^{\kern 0.5pt 
a}_{\bar a\kern 0.2pt p}\,.
\mytag{8.17}
$$
Due to the same identity \mythetag{8.17} the Lagrangian density $L$ 
itself is a real-valued scalar field in $M$.\par
     The formula \mythetag{8.16} is not new. There are some other papers, 
where the energy-momentum tensor for a spinor field is calculated, e\.\,g\.
\mycite{7}. The formula \mythetag{8.16} is quite similar to the
formula \thetag{3.15} in \mycite{7}. These formulas should 
coincide upon passing from SGS units, which are used in present
paper, to special units, where $\hbar=c=1$ instead of 
\mythetag{1.3}.\par
\enddocument
\end